\documentclass[review]{elsarticle}
\RequirePackage{amssymb}
\usepackage{bm}
\RequirePackage{amsmath}
\usepackage{amsthm}
\usepackage{color}
\usepackage{lineno,hyperref}
\modulolinenumbers[5]
\usepackage{geometry}
\usepackage{lineno,hyperref}
\modulolinenumbers[5]

\journal{~~}

\begin{document}
\newtheorem{theorem}{Theorem}[section]
\newtheorem{definition}{Definition}[section]
\newtheorem{lemma}{Lemma}[section]
\newtheorem{problem}{Problem}[section]
\newtheorem{claim}{Claim}[section]
\newtheorem{corollary}{Corollary}[section]
\newtheorem{proposition}{Proposition}[section]
\newtheorem{conjecture}{Conjecture}[section]
 \abovedisplayskip 6pt plus 2pt minus 2pt \belowdisplayskip 6pt
 plus 2pt minus 2pt
\def\vsp{\vspace{1mm}}
\def\th#1{\vspace{1mm}\noindent{\bf #1}\quad}
\def\proof{\vspace{1mm}\noindent{\it Proof}\quad}
\def\no{\nonumber}
\newenvironment{prof}[1][Proof]{\noindent\textit{#1}\quad }
{\hfill $\Box$\vspace{0.7mm}}
\begin{frontmatter}

\title{Anti-{Ramsey} Numbers of Expansions of  Doubly  Edge-critical Graphs in Uniform Hypergraphs}
\author[ams,ucas]{Tong Li\fnref{NSF3}}
\ead{litong@amss.ac.cn}
\fntext[NSF3]{Research  supported by  National Natural Science Foundation of China (Grant No. 12301459)}

\author[nuaa,klmmhpcav]{Yucong Tang\fnref{NSF1}}
\ead{tangyucong@nuaa.edu.cn}
\fntext[NSF1]{Research  supported by  National Natural Science Foundation of China (Grant No. 11901292)}

\author[ams,ucas]{Guiying Yan\corref{cor1}\fnref{NSF2}}
\ead{yangy@amss.ac.cn}
\fntext[NSF2]{Research   supported by  National Natural Science Foundation of China (Grant No. 11631014)}

\cortext[cor1]{Corresponding author}

\address[ams]{Academy of Mathematics and Systems Science, Chinese Academy of
Sciences,
 \\
 Beijing $100190$, P. R. China}

\address[ucas]
{School of Mathematical Sciences, University of Chinese Academy of Sciences, Beijing $100049$, P. R. China}

\address[nuaa]
{School of Mathematics, Nanjing University of Aeronautics and Astronautics, Nanjing $211106$, P. R. China}

\address[klmmhpcav]
{Key Laboratory of Mathematical Modelling and High Performance Computing of Air Vehicles (NUAA), MIIT,\\ Nanjing $211106$, P. R. China}

\begin{abstract}
 For an $r$-graph $H$,  the anti-Ramsey number ${\rm ar}(n,r,H)$
is the minimum number $c$ of colors such that for any edge-coloring  of the complete $r$-graph on $n$ vertices with at least $c$ colors, there is a copy of $H$ whose edges have distinct colors. A 2-graph $F$ is  doubly edge-$p$-critical  if the chromatic number $\chi(F - e)\geq p$ for every edge $e$ in $F$ and there exist two edges $e_1,e_2$ in $F$ such that $\chi(F -e_1- e_2)=p-1$. The anti-Ramsey numbers of  doubly edge-$p$-critical  2-graphs were determined by Jiang and Pikhurko \cite{Jiang&Pikhurko2009}, which generalized the anti-Ramsey numbers of cliques determined by Erd\H{o}s, Simonovits and S\'{o}s \cite{Erdos&Simonovits&Sos1975}. In general, few exact values of anti-Ramsey numbers of $r$-graphs are known for $r\geq 3$. Given a 2-graph $F$, the expansion $F^{(r)}$ of $F$ is an $r$-graph on $|V(F)|+(r-2)|F|$ vertices obtained from $F$ by adding $r-2$ new vertices to each edge of $F$. In this paper,  we determine the exact value of  ${\rm ar}(n,r,F^{(r)})$  for any doubly edge-$p$-critical 2-graph $F$ with  $p>r\geq 3$ and sufficiently large $n$.
\end{abstract}


\end{frontmatter}


\section{Introduction}

An $r$-\emph{graph} (or  $r$-\emph{uniform hypergraph})  $H$ consists of a vertex set $V$ and an edge set  with exactly $r$ vertices in each edge. For a finite set $V$ and a positive integer $r$, let ${V\choose r}$ be the set of all $r$-element subsets of $V$. We identify an $r$-graph $H$ with its edge set, i.e., $H\subseteq {{V}\choose r}$. Let $[n]=\{1,\ldots, n\}$. The $r$-graph $K_n^r={{[n]}\choose r}$ is called  complete $r$-graph on $n$ vertices. Fix an edge-coloring
of an $r$-graph $H$, a subgraph $G\subseteq H$ is \emph{rainbow} if  $G$ contains no two edges of the same color. The anti-Ramsey number of an $r$-graph $F$, denoted by ${\rm ar}(n,r,F)$,
is the minimum number $c$ of colors such that   $K_n^r$ equipped with any  edge-coloring   with at least $c$ colors contains a rainbow subgraph isomorphic to $F$.

The Tur\'{a}n number of a family of $r$-graphs $\mathcal{F}$, denoted by ${\rm ex}(n,r,\mathcal{F})$, is the maximum number of edges in an  $n$-vertex $r$-graph that does not contain any $r$-graph in $\mathcal{F}$ as a subgraph. For a given $r$-graph $F$, there is a natural lower bound of ${\rm ar}(n,r,F)$ in terms of  Tur\'{a}n number:
$${\rm ar}(n,r,F)\geq {\rm ex}(n,r,\{F-e:e\in F\}+2.$$
This trivial lower bound is easily obtained by coloring a rainbow Tur\'{a}n extremal $r$-graph for $\{F-e:e\in F\}$ in $K_n^r$, and the remaining edges with an additional color.

In 1973, Erd\H{o}s,  Simonovits, and S$\acute{\rm o}$s  \cite{Erdos&Simonovits&Sos1975} began the study of anti-Ramsey number, and proved that ${\rm ar}(n,2,K_p^2)={\rm ex}(n,2,K_{p-1}^2)+2$ for $p\geq  3$ and sufficiently large $n$. This result was extended by Montellano-Ballesteros and Neumann-Lara
\cite{Montellano-Ballesteros&Neumann-Lara2002} to all values of $n$ and $p$ with $n > p \geq 3$. For $r=2,$ the field has
been studied extensively and please see the survey \cite{Fujita&Magnant&Ozeki2010} for more details.

Given a 2-graph $F$, the \emph{chromatic number}  of  $F$, denoted by $\chi(F)$, is the smallest integer $s$ such that there is a vertex partition of $V(F)$ into $V_1\cup\cdots \cup V_{s}$ satisfying that the induced subgraph on each $V_i$ has no edge, for $i\in[s]$.    $F$ is \emph{edge-critical} if there exists an edge $e$ in $F$ such that $\chi(F) > \chi(F - e)$. $F$ is \emph{doubly edge-$p$-critical} if $\chi(F - e)\geq p$ for every edge $e$ in $F$ and there exist two edges $e_1,e_2$ in $F$ such that $\chi(F -e_1- e_2)=p-1$. Let $\mathcal{P}=V_1\cup\cdots \cup V_{p-1}$ be a $(p-1)$-partition of $V(F)$. For $i\in [p-1]$, denote $e(V_i)$  the number of edges contained in $V_i$. We call $\bm {a}_F(\mathcal{P})=(e(V_1),\ldots,e(V_{p-1}))$ the \emph{index vector} of $F$ with respect to  $\mathcal{P}$. Let $\bm{\mathcal{P}}(F)=\{\mathcal{P}: \mathcal{P} \mbox{ is a $(p-1)$-partition of $V(F)$}, \|\bm{a}_F (\mathcal{P}) \|_{\infty}=1\}$, where $\|\cdot\|_{\infty}$ is the $\infty$-norm of a  vector which equals to the maximum absolute value over all components.    Given an integer $\ell$ with $1\leq \ell\leq p-1$, a doubly edge-$p$-critical 2-graph  $F$ is  in \emph{class} $\ell$ if $\min_{\mathcal{P}\in \bm{\mathcal{P}}(F)}\|\bm{a}_F(\mathcal{P})\|_1=\ell$, where $\|\cdot\|_1$ is the $1$-norm of a  vector which equals to the sum of the absolute values of all components. When $\bm{\mathcal{P}}(F)=\emptyset$, we say $F$ is in class $p$. By the definition of doubly edge-$p$-critical graphs, if $F$ is in class $\ell$, then $2\leq \ell\leq p$. In fact, $\ell=1$ implies that there is an edge  $e\in F$ such that $\chi(F-e)=p-1$.
Let $T_p(n,r)$ be the $r$-graph with vertex set $[n]$ obtained by partitioning $[n]$ into $p-1$ parts $V_1\cup\cdots\cup V_{p-1}$ with $|V_1|\leq |V_2|\leq \cdots \leq |V_{p-1}|\leq |V_1|+1$, and the edge set of $T_p(n,r)$ consists of all the $r$-sets of vertices that  intersect every part in at most one vertex.  Let $t_p(n,r)=|T_p(n,r)|$.
Jiang and Pikhurko \cite{Jiang&Pikhurko2009} determined the exact anti-Ramsey numbers of all doubly edge-$p$-critical 2-graphs for sufficiently large $n$.

\begin{theorem}[\cite{Jiang&Pikhurko2009}]\label{2-graphantiedcritcal}
Let $n,p,\ell$ be integers such that $p>2,$ $2\leq \ell \leq p$. Given a doubly edge-$p$-critical 2-graph $F$, if $F$ is in class $\ell$, then for sufficiently large $n$,
${\rm ar}(n,2,F)=t_p(n,2)+\ell.$
\end{theorem}
Notice that $K_p^2$ is doubly edge-$p$-critical in class 2, so the above theorem is an extension of the result of Erd\H{o}s,  Simonovits, and S$\acute{\rm o}$s  \cite{Erdos&Simonovits&Sos1975}.

 For  $r\geq3$,   \"{O}zkahya and Young \cite{Ozkahya&Young2013} initiated the study of anti-Ramsey number for matchings, where a matching $M_k$ is a collection of $k$ pairwise disjoint edges.  They conjectured that
\begin{align*}{\rm ar}(n,r,M_k)=\left\{\begin{array}{ll}
{\rm ex}(n, r,M_{k-1}) + 2, & if\ k\leq c_r,n=kr,\\{\rm ex}(n, r,M_{k-1}) + r+1 , & if\ k\geq c_r,n=kr,\\
{\rm ex}(n, r,M_{k-1}) + 2, & if\ n>kr,
\end{array}\right.\end{align*}
 where $c_r$ is a constant depending on $r$.  They proved that the conjecture holds for $k=2,3$ and  sufficiently large $n$. Later,  Frankl and Kupavskii \cite{Frankl&Kupavskii2019} proved that the conjecture is true for $ n \geq rk + (r - 1)(k - 1)$ and $k \geq  3$.

Given a 2-graph $F$, the expansion $F^{(r)}$ of $F$ is an $r$-graph on   $|V(F)|+(r-2)|F|$ vertices  obtained from $F$ by adding $r-2$ new vertices to each edge in $F$.
Let $P_k$ and $C_k$ be the path and cycle with $k$ edges in 2-graphs, respectively. Recently, we \cite{Tang&Li&Yan2022,Tang&Li&Yan2022+} determined the exact values of ${\rm ar}(n,r,P_k^{(r)})$ and ${\rm ar}(n,r,C_k^{(r)})$ for all $r,k$ with $r\geq 3$, $k\geq 3$ and sufficiently large $n$, which are extensions of several results of Gu, Li and Shi \cite{Gu&Li&Shi2020}.

The motivation of this paper is to extend Theorem \ref{2-graphantiedcritcal} to $r$-graphs.  We write $x\ll y$ to mean that for any $y\in (0,1]$ there exists an $x_0\in (0,1)$, where $x_0$ is often regarded as a positive function of $y$, such that for all $x\leq x_0$ the subsequent statements hold. Hierarchies with more constants are defined in a similar way and are to be read from the right to the left. The followings are our main results.

\begin{theorem}\label{3uniform}
Let $n,p,\ell$ be integers such that $p>3,$ $2\leq \ell \leq p$ and $1/n\ll p$. Given a doubly edge-$p$-critical 2-graph $F$, if $F$ is in class $\ell$, then
${\rm ar}(n,3,F^{(3)})=t_p(n,3)+\ell.$
\end{theorem}

\begin{theorem}\label{rgeq4uniform}
Let $n,r,p,\ell$ be integers such that $p>r\geq 4,$ $2\leq \ell \leq p$ and $1/n\ll r,p$. Given a doubly edge-$p$-critical 2-graph $F$, if $F$ is in class $\ell$, then
${\rm ar}(n,r,F^{(r)})=t_p(n,r)+2.$
\end{theorem}

The idea in the proofs of Theorems \ref{3uniform} and \ref{rgeq4uniform}  is to first use the stability results for expansion of edge-critical graphs   and then modify the method of Pikhurko \cite{Pikhurko} to anti-Ramsey numbers.
The rest of this paper is organized as follows. Firstly,  we give necessary definitions, conceptions, and key lemmas for the proofs. Then  we prove Theorems \ref{3uniform} and   \ref{rgeq4uniform}.

\section{Preliminaries}

Fix $r$-graphs $G$ and $F$.  For any vertex set $U\subseteq V(G)$, define $d_G(U)=|\{e\in G: U\subseteq e\}|$, that is,  $d_G(U)$ is the number of edges in $G$ that contain $U$.
An edge-coloring of $G$ is a mapping $C: G\rightarrow [|G|]$. Given any subgraph $H\subseteq G$, let $C(H)=\{C(e):e\in H\}$ be the  set of colors appearing in $H$. Thus, $H$ is rainbow means that $C(e)\neq C(e')$ for any two edges $e,e'$ in $H$. If $H$ is isomorphic to $F$, we say $H$ is an $F$ for short sometimes.   $G$ is $F$-\emph{free} if $G$ contains no subgraph isomorphic to $F$. Fix   an edge-coloring $C$ of $G$,   $G$ is \emph{rainbow $F$-free} if $G$ has no rainbow subgraph  isomorphic to $F$.

\begin{lemma}\label{config}
Let $p,\ell$ be integers such that $p\geq 4$ and $3\leq\ell\leq p$. Let $F$ be a doubly edge-$p$-critical 2-graph in class $\ell$. Then, there is a vertex partition $\mathcal{P}=V_1\cup\cdots\cup V_{p-1}$ of $V(F)$ such that $\bm{a}_F(\mathcal{P})=(2,0,\ldots,0)$.
\end{lemma}

\begin{prof}By the definition of doubly edge-$p$-critical graph, there is an edge $e_1\in F$ such that $F-e_1$ is edge-critical and $\chi(F-e_1)=p$. Moreover, there is another edge $e_2\in F$ such that $\chi(F-e_1-e_2)=p-1$. It follows that there is a vertex partition $\mathcal{P}=V_1\cup\cdots\cup V_{p-1}$ of $F$ such that $\bm{a}_{F-e_1}(\mathcal{P})=(1,0,\ldots,0)$. As $\ell\geq 3$, $e_1$ must be contained  in $V_1$, which means that $\bm{a}_F(\mathcal{P})=(2,0,\ldots,0)$.
\end{prof}

Given an edge-colored  $K_n^r$ and a family of $r$-graphs $\mathcal{H}$. A \emph{rainbow collection} $\bm{\mathcal{H}}$ of $\mathcal{H}$ in $K_n^r$ is a collection of rainbow subgraphs such that

\noindent {\rm (i)} For any $H_1,H_2\in \bm{\mathcal{H}}$, $C(H_1)\cap C(H_2)=\emptyset$;

 \noindent {\rm (ii)} For any $H\in \bm{\mathcal{H}}$, $H$ is a rainbow copy of some $r$-graph in $\mathcal{H}$.

\begin{lemma}\label{numofrbcopy}
Let $F$ be a $2$-graph and $n,r$ be integers such that $r\geq 3$ and $1/n\ll r$. Given any edge-coloring of $K_n^r$, if $K_n^r$ is rainbow $F^{(r)}$-free, then for any $H\in \mathcal{F}^-=\{F-e:e\in F\}$ and any rainbow collection $\bm{\mathcal{H}}$ of $H^{(r)}$, we have $|\bm{\mathcal{H}}|\leq {n\choose 2}$.
\end{lemma}

\begin{prof} By the definition of $\mathcal{F}^-$, for any $H^{(r)}$ with $H\in \mathcal{F}^-$, there exists a pair of vertices $\{v_1,v_2\}$ in $V(K_n^r)$, such that $H^{(r)}\cup \{v_1,v_2,w_1,\ldots,w_{r-2}\}$ is an $F^{(r)}$, where $w_i\notin V(H^{(r)})$ for $i\in [r-2]$. We call such a pair a \emph{terminal pair} of $H^{(r)}$. Suppose $|\bm{\mathcal{H}}|> {n\choose 2}$, then the total number of terminal pairs is greater than ${n\choose 2}$. Thus, there exist $H_1^{(r)},H_2^{(r)}\in \bm{\mathcal{H}}$ sharing a terminal pair, say $\{u_1,u_2\}$. Choose $e=\{u_1,u_2,f\}$, where $f\in {{V\backslash(V(H_1^{(r)})\cup V(H_2^{(r)}))}\choose {r-2}}$. Since $C(H_1^{(r)})\cap C(H_2^{(r)})=\emptyset,$ at least one of $\{H_1^{(r)},e\}$ and $\{H_2^{(r)},e\}$ is rainbow. Clearly, $\{H_1^{(r)},e\}$ and $\{H_2^{(r)},e\}$  are copies of $F^{(r)}$, which contradicts to that $K_n^r$ is rainbow $F^{(r)}$-free.
\end{prof}

Given any positive constant $\varepsilon$ and two $r$-graphs $G_1$ and $G_2$ of  order $n$, $G_1$ is  $\varepsilon$-\emph{close} to $G_2$ if we can add or remove at most $\varepsilon n^r$ edges
from   $G_1$ to make it isomorphic to $G_2$.

The following stability result of expansions is a direct corollary of Lemma 4 in \cite{Pikhurko}.

 \begin{lemma}[\cite{Pikhurko}]\label{turanstability} Fix integers $p > r\geq 2$ and 2-graph $F$ with $\chi(F)=p$. For  $1/n\ll \eta\ll\varepsilon,p,r$, if  $G$ is an $n$-vertex $F^{(r)}$-free $r$-graph  with  $|G|\geq |T_p (n,r)|-\eta n^r$, then $G$ is $\varepsilon$-close to $T_p(n,r)$.
\end{lemma}

We use the above lemma to prove the following anti-Ramsey-type stability lemma.

\begin{lemma}\label{stability}
Let $n,r,p$ be integers such that $p>r\geq 3$ and $F$ be a  doubly edge-$p$-critical 2-graph. Let $\varepsilon, \delta$ be any positive constants with $1/n\ll \delta \ll \varepsilon$. Given an edge-colored $K_n^r$, if $K_n^r$ is rainbow $F^{(r)}$-free, then for any rainbow subgraph $G$ of  $K_n^r$ with $|G|\geq t_p(n,r)-\delta n^r$, $G$ is $\varepsilon$-close to $T_p(n,r)$.
\end{lemma}

\begin{prof}By the definition of doubly edge-$p$-critical graphs, there is an edge $e$ in $F$ such that $\chi(F-e)=p$ and $F-e$ is edge-critical. Let $\bm{\mathcal{H}}$ be a maximal edge-disjoint collection of $(F-e)^{(r)}$  in $G$, i.e., (i) any two $r$-graphs in $\bm{\mathcal{H}}$ are edge-disjoint; (ii) there is no copy of $(F-e)^{(r)}$ in $G$ after deleting all the edges in $\bm{\mathcal{H}}$. Since $G$ is rainbow, $\bm{\mathcal{H}}$ is rainbow collections of $(F-e)^{(r)}$ in $K_n^r$. By Lemma \ref{numofrbcopy}, $|\bm{\mathcal{H}}|\leq {n\choose 2}$. Let $G'$ be the subgraph of $G$ by deleting from $G$ all the edges in $\bm{\mathcal{H}}$. Thus, $G'$ is $(F-e)^{(r)}$-free and $|G'|\geq |G|-|F|{n\choose 2}\geq t_p(n,r)-2\delta n^r$. Therefore, by applying Lemma \ref{turanstability} with $2\delta=\eta$, $G'$ is $\frac{\varepsilon}{2}$-close to $T_p(n,r)$, which follows that $G$ is $\varepsilon$-close to $T_p(n,r)$.
\end{prof}

\section{3-uniform hypergraphs}
Through out the proofs, we use the following definitions and constants. Let $\tau = |V(F^{(r)})|$ and $v(F)=|V(F)|$. A pair of vertices $\{u,v\}$ is called  {\it big} in an $r$-graph $G$ if $d_G(u,v)\geq \tau{n\choose {r-3}}+|F|$, and \emph{small} otherwise. $F$ is called the \emph{skeleton} of $F^{(r)}$. It is easy to see that if $G$ is rainbow and there is a skeleton of $F^{(r)}$ such that all the vertex pairs in the skeleton  are big in $G$, then one can greedily extend the skeleton to a rainbow $F^{(r)}$.

\begin{prof}[Proof of Theorem \ref{3uniform}]  First, we prove ${\rm ar}(n,3,F^{(3)})> t_p(n,3)+\ell-1$ by giving an extremal coloring. Color a  $T_p(n,3)$ in $K_n^3$ into a rainbow $T_p(n,3)$. Let $V_1\cup\cdots\cup V_{p-1}$ be the vertex partition corresponding to $T_p(n,3)$. Then, for $i\in[\ell-2]$, let $E_i=\{e\in K_n^3: |e\cap V_i|\geq 2\}$ and color each $E_i$ with distinct new colors and the remaining edges with a new color. Thus, the total number of colors is $t_p(n,3)+\ell-2+1$. Suppose that there is a rainbow $F^{(3)}$ in $K_n^3$. Let $U_i=V(F)\cap V_i$, $i\in[p-1]$ and $\mathcal{P}=\cup_{i\in[p-1]}U_i$. Notice that for $i\in [\ell-2]$, $E_i\cap T_p(n,3)=\emptyset$ and the edges in $E_i$ have the same color, so the number of edges of $F$ contained in $U_i$ is at most one. As $\chi(F)=p$, $\bm{a}_F(\mathcal{P})\neq (0,\ldots,0)$. Hence, $F$ is in class $\ell'$, where $\ell'\leq \ell-1$, a contradiction.

Now, we prove the upper bound. Let $\varepsilon_i$, $i\in [7]$ be positive constants with $1/n\ll\varepsilon_7\ll\cdots \ll \varepsilon_1\ll r,p$. Let $K_n^3$ be edge-colored with $t_p(n,3)+\ell$ colors and rainbow $F^{(3)}$-free. Let $G$ be the subgraph of $K_n^3$ obtained by taking one edge from each color class. Thus $|G|=t_p(n,3)+\ell$. By Lemma  \ref{stability}, $G$ is $\varepsilon_7$-close to $T_p(n,3)$.

For any vertex partition $\mathcal{P}=V_1\cup \cdots \cup V_{p-1}$ of $V(G)$, define $f_G(\mathcal{P})=\sum_{e\in G}|\{i\in[p-1]:|e\cap V_i|\neq \emptyset\}|$ and let $\mathcal{P}^*$ be the vertex partition that maximize $f_G(\mathcal{P})$. That is $$f_G(\mathcal{P}^*)=\max\{f_G(\mathcal{P}):{\mathcal{P}\mbox{ is a $(p-1)$-partition of $V(G)$}}\}.$$

Set $\mathcal{P}^*=V_1\cup\cdots\cup V_{p-1}$ and $T=\{e\in {{V(G)}\choose 3}: |e\cap V_i|\leq 1,i\in[p-1]\}$. Then
\begin{align}\label{1}f_G(\mathcal{P}^*)\leq 3|G|-|G\backslash T|.
\end{align}
 Let $\mathcal{P}'$ be the vertex partition corresponding to $T_p(n,3)$. We have
\begin{align}\label{2}f_G(\mathcal{P}^*)\geq f_G(\mathcal{P}') \geq 3|G\cap T_p(n,3)|\geq 3(|G|-\varepsilon_7 n^3),
\end{align}
where the last inequality is due to that $G$ is $\varepsilon_7$-close to $T_p(n,3)$. So, by  (\ref{1}) and (\ref{2}), $|G\backslash T|\leq 3\varepsilon_7 n^3$. Therefore, $$|T|\geq |T\cap G|=|G|-|G\backslash T|\geq t_p(n,3)+\ell-3\varepsilon_7 n^3.$$
This lower bound on $|T|$ can be shown to imply (see Claim 1 in \cite[Proof of Theorem 5]{Mubayi1}) that $|V_i|\geq \frac{n}{2(p-1)}$  for $i\in[p-1]$.

Let us call edges in $G\backslash T$ \emph{non-crossing}. For a non-crossing edge $e$, call $\{\{u,v\}:u,v\in e\cap V_i \mbox{ for some } i\in[p-1]\}$ the set of \emph{non-crossing parts} of $e$. Next, the number of non-crossing edges is lower bounded as follows.
\begin{align*}
|G\backslash T|&=|G|-|G\cap T |\geq t_p(n,3)+\ell-|G\cap T|\\
&\geq |T|+\ell-|G\cap T|\\
&=|T\backslash G|+\ell.
\end{align*}
It follows that the number of non-crossing edges in $G$ is at least $\ell$. And we can bound the number of small pairs.

\begin{claim}\label{3uniformclaim1}The number of small pairs  is at most $\varepsilon_6 n^2$.
\end{claim}

\begin{prof}Let $x$ be the number of small pairs in $G$. Since each small pair corresponds to at least $\big[(p-3)\frac{n}{2(p-1)}-(\tau+|F|)\big]$ edges in $T\backslash G$, and each edge  in $T\backslash G$ contains at most three small pairs, we have
\begin{align*}x\frac{\big[(p-3)\frac{n}{2(p-1)}-(\tau+|F|)\big]}{3}\leq |T\backslash G|\leq |G\backslash T|-\ell\leq 2\varepsilon_7n^3,
\end{align*}
which implies $x\leq \varepsilon_6 n^2$.
\end{prof}

Let $\mathcal{Q}=U_1\cup\cdots U_{p-1}$ be the vertex partition of $V(F)$ with $\|\bm{a}_F(\mathcal{Q})\|_1=\ell$ and $q_i=|U_i|$ for $i\in[p-1]$. If $\ell\geq 3$, then by Lemma \ref{config}, there is a vertex partition $\mathcal{Q'}=U_1'\cup\cdots\cup U_{p-1}'$ of $V(F)$ such that $\bm{a}_F(\mathcal{Q'})=(2,0,\ldots,0)$. Let $q'_i=|U_i'|$ for $i\in[p-1]$. For any set $\{e_1,\ldots,e_{\ell}\} $    of $\ell$ non-crossing edges with $\{u_{2i-1},u_{2i}\}$ being  a   non-crossing part of $e_i$ for $i\in [\ell]$, define $\Omega$ to be the Cartesian product of vertex sets as follow, i.e.,
$$\Omega={{V_1'}\choose {v(F)+2}}\times {{V_2'}\choose {v(F)+2}}\times\cdots \times {{V_{p-1}'}\choose {v(F)+2}},$$
where $V_i'=V_i\backslash (\cup_{i\in [\ell]} e_i)$. The element in $\Omega$ can be regarded as choices of vertices used to form the skeleton of $F^{(3)}$. Let $\Omega'$ be  the set of all choices $(X_1,\ldots,X_{p-1})$ in $\Omega$ such that ${{\cup_{i\in[p-1]}X_i}\choose 2}$ contains at least one small pair. For  $u\in\{u_i:i\in[2\ell]\}$, let $\Omega_u$ be the  subset of all choices $(X_1,\ldots,X_{p-1})$ in $\Omega$ such that there is a small pair containing $u$ and a vertex in $\cup_{i\in[p-1]}X_i$.

\begin{claim}\label{3uniformclaimnew} $\Omega=\Omega'\cup (\cup_{i\in[2\ell]}\Omega_{u_i})$.
\end{claim}
\begin{prof}
Let   $\mathcal{X}=(X_1,\ldots,X_{p-1})\in \Omega$, where $X_i=\{x_{i1},\ldots,x_{i,v(F)},y_{i1},y_{i2}\}$ for $i\in [p-1]$. We split the proof into several cases depending on $\ell$ and the positions of $\{u_{2i-1},u_{2i}\}_{i\in [\ell]}$.

\noindent{\bf Case 1.} \quad $\{u_{2i-1},u_{2i}\}$, $i\in [\ell]$,  are in different vertex classes.

Without loss of generality, assume $\{u_{2i-1},u_{2i}\} \subseteq V_i$ for $i\in[\ell]$.  Notice that in this case, $\ell\leq p-1.$

\noindent{\bf Subcase 1.1}\quad $\ell=2$.

 Consider the two edges $e_1'=\{x_{11},x_{12},y_{11}\}$, $e_2'=\{x_{21},x_{22},y_{21}\}$ in $K_n^3$. If $C(e_1')=C(e_2)$ and $C(e_2')=C(e_1)$, then $C(e_1')\neq C(e_2')$. Thus, there is a small pair in ${{\cup_{i\in[p-1]}X_i}\choose 2}$, that is, $\mathcal{X}\in \Omega'$. Otherwise, let $X_i'=\{x_{i1},\ldots,x_{i,q_i}\}$ for $i\in [p-1]$. Then, $\cup_{i\in[p-1]}X_i'$  could form the skeleton of $F^{(3)}$, and then be expanded to a rainbow copy of $F^{(3)}$ containing $e_1'$ and $e_2'$, a contradiction. If $C(e_1')\neq C(e_2)$, then  there is a small pair consisting of a vertex in $\{u_3,u_4\}$ and another vertex in $\cup_{i\in[p-1]}X_i$, that is, $\mathcal{X}\in \Omega_{u_3}\cup \Omega_{u_4}$. Otherwise, let $X_i'=\{x_{i1},\ldots,x_{i,q_i}\}$ for $i\in [p-1]\backslash \{2\}$ and $X_2'=\{x_{21},\ldots,x_{2,q_2-2}\}$. Then, $\cup_{i\in[p-1]}X_i'\cup\{u_3,u_4\}$ could form the skeleton of $F^{(3)}$, and then be expanded to a rainbow copy of $F^{(3)}$ containing $e_1'$ and $e_2$, a contradiction. The case $C(e_2')\neq C(e_1)$ is similar.

Thus, $\Omega=\Omega'\cup (\cup_{i\in[2\ell]}\Omega_{u_i})$.

\noindent{\bf Subcase 1.2}\quad $3\leq \ell\leq p-1$ and the two edges of $F$ in $U_1'$ intersect.

 Consider the $2\ell$ edges $e_i'=\{x_{i1},x_{i2},y_{i1}\}$ and $e_i''=\{u_{2i},x_{i1},y_{i2}\}$, $i\in[\ell]$, in $K_n^3$. If $C(e_i')=C(e_i)$ holds for $i\in [\ell]$, then $e_1',\ldots,e_{\ell}'$ is rainbow. Thus, $\mathcal{X}\in \Omega'$. Otherwise, let $X_i'=\{x_{i1},\ldots,x_{i,q_i}\}$. Then, $\cup_{i\in[p-1]}X_i'$  could form the skeleton of $F^{(3)}$, and then be expanded to a rainbow copy of $F^{(3)}$ containing $\{e_1',\ldots,e_{\ell}'\}$. If $C(e_1')\neq C(e_1)$, then at least one of the inequalities  $C(e_1'')\neq C(e_1)$ and $C(e_1')\neq C(e_1'')$ holds. If $C(e_1'')\neq C(e_1)$, then $\mathcal{X}\in \Omega_{u_1}\cup \Omega_{u_2}$. Otherwise,  let $X_i'=\{x_{i1},\ldots,x_{i,q'_i}\}$ for $i\in [p-1]\backslash \{1\}$ and $X_1'=\{x_{11},\ldots,x_{1,q'_1-2}\}$. Then,  $\cup_{i\in[p-1]}X_i'\cup\{u_1,u_2\}$ could form the skeleton of $F^{(3)}$, and then be expanded to a rainbow copy of $F^{(3)}$ containing $e_1''$ and $e_1$. Similarly, if $C(e_1')\neq C(e_1'')$, then $\mathcal{X}\in \Omega_{u_2}$.

Thus, $\Omega=\Omega'\cup (\cup_{i\in[2\ell]}\Omega_{u_i})$.

\noindent{\bf Subcase 1.3}\quad $3\leq \ell\leq p-1$ and the two edges of $F$ in $U_1'$ are disjoint.

 Consider the $\ell$ edges $e_i'=\{x_{i1},x_{i2},y_{i1}\}$, $i\in[\ell]$, in $K_n^3$. If $C(e_i')=C(e_i)$ holds for $i\in [\ell]$, then $e_1',\ldots,e_{\ell}'$ is rainbow. Thus,  $\mathcal{X}\in \Omega'$. Otherwise, let $X_i'=\{x_{i1},\ldots,x_{i,q_i}\}$. Then, $\cup_{i\in[p-1]}X_i'$  could form the skeleton of $F^{(3)}$, and then be expanded to a rainbow copy of $F^{(3)}$ containing $\{e_1',\ldots,e_{\ell}'\}$.  If $C(e_1')\neq C(e_1)$, then  $\mathcal{X}\in \Omega_{u_1}\cup \Omega_{u_2}$, which is similar to the Subcase 1.2.

Thus, $\Omega=\Omega'\cup (\cup_{i\in[2\ell]}\Omega_{u_i})$.

\noindent{\bf Case 2.} \quad There are two non-crossing parts contained in a common vertex class.

Without loss of generality, assume $\{u_1,\ldots,u_{4}\} \subseteq V_1$ and $u_1\neq u_3$, $u_1\neq u_4$, $u_2\neq u_3$. Let $V_i'=V_i\backslash ( e_1\cup e_2)$.

\noindent{\bf Subcase 2.1} \quad $\ell= 2$.

  Consider $e=\{x_{21},x_{22},y_{21}\}$. At least one of $\{e,e_1\}$ and  $\{e,e_2\}$ is rainbow. Thus, $\mathcal{X}\in\Omega'\cup (\cup_{i\in[2\ell]}\Omega_{u_i})$. Otherwise, let $X_i'=\{x_{i1},\ldots,x_{i,q_i}\}$ for $i\in [p-1]\backslash \{1\}$ and $X_1'=\{x_{11},\ldots,x_{1,q_1-2}\}$. It follows  that   both  $\cup_{i\in[p-1]}X_i'\cup \{u_1,u_2\}$ and $\cup_{i\in[p-1]}X_i'\cup \{u_3,u_4\}$ could form the skeleton of $F^{(3)}$. Thus at least one of them can be expanded to a rainbow copy of $F^{(3)}$ containing either $\{e, e_1\}$ or $\{e,e_2\},$ a contradiction.

\noindent{\bf Subcase 2.2} \quad $\ell\geq 3$ and the two edges of $F$ in $U_1'$ intersect.

Consider $e=\{u_1,u_3,y_{11}\}$. Then at least one of $\{e,e_1\}$ and  $\{e,e_2\}$ is rainbow. Thus, $\mathcal{X}\in\Omega'\cup (\cup_{i\in[2\ell]}\Omega_{u_i})$. Otherwise, let $X_i'=\{x_{i1},\ldots,x_{i,q'_i}\}$ for $i\in [p-1]\backslash \{1\}$ and $X_1'=\{x_{11},\ldots,x_{1,q'_1-3}\}$. It follows  that both  $\cup_{i\in[p-1]}X_i'\cup \{u_1,u_2,u_3\}$ and $\cup_{i\in[p-1]}X_i'\cup \{u_1,u_3,u_4\}$ could form the skeleton of $F^{(3)}$. Thus at least one of them can be expanded to a rainbow copy of $F^{(3)}$ containing either $\{e, e_1\}$ or $\{e,e_2\},$ a contradiction.

\noindent{\bf Subcase 2.3}\quad $\ell\geq 3$ and the two edges of $F$ in $U_1'$ are disjoint.

  Consider $e=\{x_{11},x_{12},y_{11}\}$. Then at least one of $\{e,e_1\}$ and  $\{e,e_2\}$ is rainbow. Thus, $\mathcal{X}\in\Omega'\cup (\cup_{i\in[2\ell]}\Omega_{u_i})$. Otherwise, let $X_i'=\{x_{i1},\ldots,x_{i,q'_i}\}$ for $i\in [p-1]\backslash \{1\}$ and $X_1'=\{x_{11},\ldots,x_{1,q'_1-2}\}$. It follows  that   both  $\cup_{i\in[p-1]}X_i'\cup\{u_1,u_2\}$ and $\cup_{i\in[p-1]}X_i'\cup\{u_3,u_4\}$ could form the skeleton of $F^{(3)}$. Thus at least one of them can be expanded to a rainbow copy of $F^{(3)}$ containing either $\{e, e_1\}$ or $\{e,e_2\},$ a contradiction.

Therefore, in all cases, we have $\Omega=\Omega'\cup (\cup_{i\in[2\ell]}\Omega_{u_i})$.
\end{prof}

For $v\in V(G)$, define $d_s(v)=|\{u\in V(G): \{u,v\} \mbox{ is small}\}|.$ Let $A=\{v\in V(G):d_s(v)\geq \varepsilon_3 n\}$.

\begin{claim}\label{3uniformclaim2} Let $\{e_1,\ldots,e_{\ell}\} $ be any set  of $\ell$ non-crossing edges with $\{u_{2i-1},u_{2i}\}$ being  a   non-crossing part of $e_i$ for $i\in [\ell]$. Then, at least one of the vertices $u_i$, $i\in[2\ell]$, is in $A$.
\end{claim}
\begin{prof}
By Claim \ref{3uniformclaimnew}, $\Omega=\Omega'\cup (\cup_{i\in[2\ell]}\Omega_{u_i})$.
Then, we bound the size of $\Omega_{u_i}$. Clearly, $$|\Omega|\geq   \bigg({{\frac{n}{3(p-1)}}\choose {v(F)+2}}\bigg)^{p-1}\geq \varepsilon_2 n^{(p-1)(v(F)+2)}.$$
Since each small pair contained in $\Omega'$ corresponds to at most $n^{(p-1)(v(F)+2)-2}$ choices in $\Omega'$, combined with Claim \ref{3uniformclaim1}, we have
$$|\Omega'|\leq \varepsilon_6 n^2\cdot n^{(p-1)(v(F)+2)-2}\leq \frac{1}{2}|\Omega|.$$
Therefore,
$$\sum_{i\in[2\ell]}|\Omega_{u_i}|\geq \frac{1}{2}|\Omega|,$$
which follows that there is some vertex $u\in\{u_i:i\in [2\ell]\}$ with $|\Omega_u|\geq \frac{1}{4\ell}|\Omega|\geq \frac{\varepsilon_2}{4\ell}n^{(p-1)(v(F)+2)}$. Hence,
$$d_s(u)\geq \frac{\varepsilon_2 n^{(p-1)(v(F)+2)}}{4\ell n^{(p-1)(v(F)+2)-1}}\geq \varepsilon_3 n.$$
\end{prof}

Now, let $G'$ be the subgraph of $G$ obtained by deleting from $G$ all the non-crossing edges whose non-crossing parts are disjoint from $A$. By Claim \ref{3uniformclaim2}, the number of deleted edges is at most $\ell-1$.

For any $v\in A$, define $N_{ncp}(v)=\{u\in V(G'): \{u,v\} \mbox{ is a non-crossing part of some non-crossing edge in } G'\}$ and $d_{ncp}(v)=|N_{ncp}(v)|$.

\begin{claim}\label{3unifclaim3}
There is a vertex $v_0$ in $A$ such that the number of non-crossing edges whose non-crossing parts containing $v_0$ is at least  $\varepsilon_4 n^2$ and $d_{ncp}(v_0)\geq \varepsilon_5 n$.
\end{claim}
\begin{prof} By the definition of $A$, the number of small pairs intersecting $A$ is at least $|A|\varepsilon_3 n/2$. Then,
$$|T\backslash G|\geq \frac{|A|\varepsilon_3 n}{2}\bigg((p-3)\frac{n}{2(p-1)}-\tau-|F|\bigg)\frac{1}{3}\geq \varepsilon_4 |A| n^2.$$
Thus,
$$|G'\backslash T|\geq |G\backslash T|-(\ell-1)\geq |T\backslash G|+\ell-(\ell-1)\geq \varepsilon_4 |A| n^2.$$
Since the non-crossing part of every non-crossing edge in $G'\backslash T$ intersects $A$,  there is a vertex $v_0\in A$ such that the number of non-crossing edges whose non-crossing parts containing $v_0$ is at least  $\varepsilon_4 n^2$. And this implies $d_{ncp}(v_0)\geq \varepsilon_4 n^2/n\geq \varepsilon_5 n$.
\end{prof}

Without loss of generality, assume $v_0\in V_1$. For $2\leq j\leq p-1$, let $Z_j=\{v\in V_j:\{v_0,v\} \mbox{ is big}\}$. For $y\in N_{ncp}(v_0)$, let $e_{v_0}(\cdot)$ be any mapping from $N_{ncp}(v_0)$ to the set of non-crossing edges containing $v_0$ such that $e_{v_0}(y)$ contains $y$. Now we finish the proof by discussion on the size of each $Z_j$, $2\leq j \leq p-1$.

\noindent{\bf Case A.} \quad $|Z_j|\geq \varepsilon_5 n$ for all $2\leq j \leq p-1$.

\noindent{\bf Subcase A.1}\quad  $\ell\geq 3$ and the two edges of $F$ in $U_1'$ intersect.

Let \begin{align*}&\Omega=\bigg\{(y,z,w,W_1,\ldots,W_{p-1})\in N_{ncp}(v_0)\times V_1\times V\times{V_1\choose {q'_1-3}}\times{Z_2\choose {q'_2}}\times\cdots \times{Z_{p-1}\choose {q'_{p-1}}}:\\
&\quad\quad \ \ \  w,z\notin e_{v_0}(y), (\{w,z\}\cup e_{v_0}(y))\cap (\cup_{i\in[p-1]}W_i)=\emptyset \bigg\},\\
&\Omega_0=\{(y,z,w,W_1,\ldots,W_{p-1})\in \Omega:
C(e_{v_0}(y))=C(\{v_0,z,w\})\},\\
& \Omega_1=\Omega\backslash\Omega_0.
\end{align*}

First, $$|\Omega|\geq\varepsilon_5 n\frac{n}{3(p-1)}\frac{n}{2}{{\frac{n}{3(p-1)}}\choose {q'_1-3}}\Pi_{i=2}^{p-1}{{\frac{\varepsilon_5n}{2}}\choose {q'_i}}>2\varepsilon_6n^{v(F)}.$$
Next, we bound $|\Omega_0|.$ Recall that $G$ is rainbow. Thus, for a fixed choice of $\{z,w\}$, there are at most two choices of $y$ such that $C(e_{v_0}(y))=C(\{v_0,z,w\})$. Therefore, $|\Omega_0|\leq 2{n\choose 2} n^{v(F)-3}\leq \frac{|\Omega|}{2}$.
As $\Omega=\Omega_0\cup\Omega_1$, we have $|\Omega_1|\geq |\Omega|/2$. Since $C(e_{v_0}(y))\neq C(\{v_0,z,w\})$ holds for every choice in $\Omega_1$, there is at least one small pair which does not contain $v_0$  in each choice in $\Omega_1$. Each such small pair is counted at most $n^{v(F)-2}$ times in $\Omega_1$. Thus, the total number of small pairs is at least
$$\frac{|\Omega_1|}{n^{v(F)-2}}\geq \frac{|\Omega|}{2n^{v(F)-2}}> \varepsilon_6 n^2 ,$$
which contradicts to Claim \ref{3uniformclaim1}.

\noindent{\bf Subcase A.2}\quad  $\ell\geq 3$ and the two edges of $F$ in $U_1'$ are disjoint.

Similar to Subcase A.1, let \begin{align*}&\Omega=\bigg\{(y,Z,w,W_1,\ldots,W_{p-1})\in N_{ncp}(v_0)\times {V_1\choose 2}\times V\times{V_1\choose {q'_1-4}}\times{Z_2\choose {q'_2}}\times\cdots \times{Z_{p-1}\choose {q'_{p-1}}}:\\
&\quad\quad \ \ \  (\{w\}\cup Z) \cap  e_{v_0}(y)=\emptyset, (\{w\}\cup Z \cup e_{v_0}(y))\cap (\cup_{i\in[p-1]}W_i)=\emptyset \bigg\},\\
&\Omega_0=\{(y,Z,w,W_1,\ldots,W_{p-1})\in \Omega:
C(e_{v_0}(y))=C(\{Z,w\})\},\\
&\Omega_1=\Omega\backslash\Omega_0.
\end{align*}

First, $$|\Omega|\geq\varepsilon_5 n{{\frac{n}{3(p-1)}}\choose {2}}\frac{n}{2}{{\frac{n}{3(p-1)}}\choose {q'_1-4}}\Pi_{i=2}^{p-1}{{\frac{\varepsilon_5n}{2}}\choose {q'_i}}>2\varepsilon_6n^{v(F)}.$$
Next, we bound $|\Omega_0|.$ Again, notice that $G$ is rainbow. Thus, for a fixed choice of $\{Z,w\}$, there are at most two choices of $y$ such that $C(e_{v_0}(y))=C(\{Z,w\})$. Therefore, $|\Omega_0|\leq 2{n\choose 3} n^{v(F)-4}\leq \frac{|\Omega|}{2}$.
As $\Omega=\Omega_0\cup\Omega_1$, we have $|\Omega_1|\geq |\Omega|/2$. Since $C(e_{v_0}(y))\neq C(\{Z,w\})$ holds for every choice in $\Omega_1$, there is at least one small pair which does not contain $v_0$  in each choice in $\Omega_1$. Each such small pair is counted at most $n^{v(F)-2}$ times in $\Omega_1$. Thus, the total number of small pairs is at least
$$\frac{|\Omega_1|}{n^{v(F)-2}}\geq \frac{|\Omega|}{2n^{v(F)-2}}> \varepsilon_6 n^2 ,$$
which contradicts to Claim \ref{3uniformclaim1}.

\noindent{\bf Subcase A.3}\quad  $\ell=2$.

In this case, we use the partition $\mathcal{Q}$ of $F$.
Similarly, let \begin{align*}&\Omega=\bigg\{ (y,Z,w,W_1,\ldots,W_{p-1})\in N_{ncp}(v_0)\times {Z_2\choose 2}\times V\times{V_1\choose {q_1-2}}\times{Z_2\choose {q_2-2}}\times {Z_3\choose {q_3}}\times\cdots \times{Z_{p-1}\choose {q_{p-1}}}:\\
&\quad\quad \ \ \  (\{w\}\cup Z) \cap  e_{v_0}(y)=\emptyset, (\{w\}\cup Z \cup e_{v_0}(y))\cap (\cup_{i\in[p-1]}W_i)=\emptyset \bigg\},\\
&\Omega_0=\{(y,Z,w,W_1,\ldots,W_{p-1})\in \Omega:
C(e_{v_0}(y))=C(\{Z,w\})\},\\
& \Omega_1=\Omega\backslash\Omega_0.
\end{align*}

First, $$|\Omega|\geq\varepsilon_5 n{{\frac{\varepsilon_5n}{2}}\choose {2}}\frac{n}{2}{{\frac{n}{3(p-1)}}\choose {q_1-2}}{{\frac{\varepsilon_5n}{2}}\choose {q_2-2}}\Pi_{i=3}^{p-1}{{\frac{\varepsilon_5n}{2}}\choose {q_i}}>2\varepsilon_6n^{v(F)}.$$
Next, we bound $|\Omega_0|.$ Again, notice that $G$ is rainbow. Thus, for a fixed choice of $\{Z,w\}$, there are at most two choices of $y$ such that $C(e_{v_0}(y))=C(\{Z,w\})$. Therefore, $|\Omega_0|\leq 2{n\choose 3} n^{v(F)-4}\leq \frac{|\Omega|}{2}$.
As $\Omega=\Omega_0\cup\Omega_1$, we have $|\Omega_1|\geq |\Omega|/2$. Since $C(e_{v_0}(y))\neq C(\{Z,w\})$ holds for every choice in $\Omega_1$, there is at least one small pair which does not contain $v_0$  in each choice in $\Omega_1$. Each such small pair is counted at most $n^{v(F)-2}$ times in $\Omega_1$. Thus, the total number of small pairs is at least
$$\frac{|\Omega_1|}{n^{v(F)-2}}\geq \frac{|\Omega|}{2n^{v(F)-2}}> \varepsilon_6 n^2 ,$$
which contradicts to Claim \ref{3uniformclaim1}.

\noindent{\bf Case B.} \quad There is some $j\in[p-1]\backslash \{1\}$ such that $|Z_j|< \varepsilon_5 n$.

Without loss of generality, assume $j=2$.
We show that if $v_0$ is moved from $V_1$ to $V_2$, then $f_G(\mathcal{P^*})$ will strictly increase, which contradicts to that $f_G(\mathcal{P^*})$ is maximum. Let $E'=\{e\in G: e\cap V_1=\{v_0\},e\cap V_2\neq \emptyset\}$, $E''=\{e\in G: v_0\in e, |e\cap V_1|\geq 2,e\cap V_2= \emptyset\}$. It suffices to show that $|E''|>|E'|$.
As $|Z_2|< \varepsilon_5 n$, we have
$$|E'|\leq \varepsilon_5 n^2+(n-\varepsilon_5 n)(\tau+|F|)\leq 2\varepsilon_5 n^2.$$
To bound $|E''|$, recall that the number of non-crossing edges with non-crossing part containing $v_0$ is at least $\varepsilon_4 n^2$ by Claim \ref{3unifclaim3}. Moreover, the number of edges containing $v_0$ and a vertex in $V_2$ is at most $\varepsilon_5 n^2+(n-\varepsilon_5 n)(\tau+|F|)\leq 2\varepsilon_5 n^2.$ Thus, $|E''|\geq \varepsilon_4 n^2-2\varepsilon_5 n^2>2\varepsilon_5 n^2=|E'|$. The proof is completed.
\end{prof}


\section{$r$-uniform hypergraphs with $r\geq 4$}
\begin{prof}[Proof of Theorem \ref{rgeq4uniform}] For the lower bound, firstly color a $T_p(n,r)$ in $K_n^r$ into a rainbow $T_p(n,r)$. Let $V_1\cup\cdots\cup V_{p-1}$ be the vertex partition corresponding to $T_p(n,r)$. Then, color the remaining edges with a new color. Thus, the total number of colors is $t_p(n,r)+1$. Suppose that there is a rainbow $F^{(r)}$ in $K_n^r$. Let $U_i=V(F)\cap V_i$, $i\in[p-1]$ and $\mathcal{P}=\cup_{i\in[p-1]}U_i$. Notice that the edges not in $ T_p(n,r)$   have the same color, so   there is at most one $U_i$ with $i\in[p-1]$ such that the number of edges of $F$ in $U_i$ is  one. As $\chi(F)=p$, $\bm{a}_F(\mathcal{P})\neq (0,\ldots,0)$. Hence, $F$ is in class $1$,  a contradiction.

Now, we prove the upper bound. Let $\varepsilon_i$, $i\in [7]$ be positive constants with $1/n\ll\varepsilon_7\ll\cdots \ll \varepsilon_1\ll r,p$. Let $K_n^r$ be edge-colored with $t_p(n,r)+2$ colors and rainbow $F^{(r)}$-free. Let $G$ be the subgraph of $K_n^r$ obtained by taking one edge from each color class. Thus $|G|=t_p(n,r)+2$. By Lemma \ref{stability}, $G$ is $\varepsilon_7$-close to $T_p(n,r)$.

For any vertex partition $\mathcal{P}=V_1\cup \cdots \cup V_{p-1}$ of $V(G)$, define $f_G(\mathcal{P})=\sum_{e\in G}|\{i\in[p-1]:|e\cap V_i|\neq \emptyset\}|$ and let $\mathcal{P}^*$ be the vertex partition that maximize $f_G(\mathcal{P})$. That is, $$f_G(\mathcal{P}^*)=\max\{f_G(\mathcal{P}):{\mathcal{P}\mbox{ is a $(p-1)$-partition of $V(G)$}}\}.$$

Set $\mathcal{P}^*=V_1\cup\cdots\cup V_{p-1}$ and $T=\{e\in {{V(G)}\choose r}: |e\cap V_i|\leq 1,i\in[p-1]\}$. Then
\begin{align}\label{3}f_G(\mathcal{P}^*)\leq r|G|-|G\backslash T|.
\end{align}
 Let $\mathcal{P}'$ be the vertex partition corresponding to $T_p(n,r)$. We have
\begin{align}\label{4}f_G(\mathcal{P}^*)\geq f_G(\mathcal{P}') \geq r|G\cap T_p(n,r)|\geq r(|G|-\varepsilon_7 n^3),
\end{align}
where the last inequality is due to that $G$ is $\varepsilon_7$-close to $T_p(n,r)$. So, by  (\ref{3}) and (\ref{4}), $|G\backslash T|\leq r \varepsilon_7 n^r$. Therefore, $$|T|\geq |T\cap G|=|G|-|G\backslash T|\geq t_p(n,r)+2-r\varepsilon_7 n^r.$$
This lower bound on $|T|$ can be shown to imply (see Claim 1 in \cite[Proof of Theorem 5]{Gu&Li&Shi2020}) $|V_i|\geq \frac{n}{2(p-1)}$  for $i\in[p-1]$.

Next, the number of non-crossing edges is lower bounded as follow.
\begin{align*}
|G\backslash T|&=|G|-|G\cap T |\geq t_p(n,r)+2-|G\cap T|\\
&\geq |T|+2-|G\cap T|\\
&=|T\backslash G|+2.
\end{align*}
It follows that the number of non-crossing edges in $G$ is at least 2. And we can bound the number of small pairs similarly as Claim \ref{3uniformclaim1} in the proof of Theorem \ref{3uniform}.

\begin{claim}\label{runiformclaim1}The number of small pairs  is at most $\varepsilon_6 n^2$.
\end{claim}

\begin{prof}Let $x$ be the number of small pairs in $G$. Since each small pair corresponds to at least ${{p-3}\choose {r-2}}\big(\frac{n}{2(p-1)}\big)^{r-2}-({n\choose {r-3}}\tau+|F|)$ edges in $T\backslash G$, and each edge in $T\backslash G$ contains at most ${r \choose 2}$ small pairs, we  have
\begin{align*}x\frac{{{p-3}\choose {r-2}}\big(\frac{n}{2(p-1)}\big)^{r-2}-({n\choose {r-3}}\tau+|F|)}{{r \choose 2}}\leq |T\backslash G|\leq |G\backslash T|-2\leq (r-1)\varepsilon_7n^r,
\end{align*}
which implies $x\leq \varepsilon_6 n^2$.
\end{prof}

Let $\mathcal{Q}=U_1\cup\cdots U_{p-1}$ be the vertex partition of $F$ with $\|\bm{a}_F(\mathcal{Q})\|_1=\ell$ and $q_i=|U_i|$ for $i\in[p-1]$. If $\ell\geq 3$, then by Lemma \ref{config}, there is a vertex partition $\mathcal{Q'}=U_1'\cup\cdots\cup U_{p-1}'$ of $F$ such that $\bm{a}_F(\mathcal{Q'})=(2,0,\ldots,0)$. Let $q'_i=|U_i'|$ for $i\in[p-1]$. For any set $\{e_1,e_2\} $    of two non-crossing edges $\{u_{2i-1},u_{2i}\}$ being  a   non-crossing part of $e_i$ for $i=1,2$, define $\Omega$ to be the Cartesian product of vertex sets as follow,
$$\Omega={{V_1'}\choose {v(F)+r-2}}\times {{V_2'}\choose {v(F)+r-2}}\times{{V_{3}'}\choose {v(F)}}\times\cdots \times {{V_{p-1}'}\choose {v(F)}},$$
where $V_i'=V_i\backslash (e_1 \cup e_2)$.  Let $\Omega'$ be  the set of all choices $(X_1,\ldots,X_{p-1})$ in $\Omega$ such that ${{\cup_{i\in[p-1]}X_i}\choose 2}$ contains at least one small pair. For  $u\in\{u_i:i\in[4]\}$, let $\Omega_u$ be the  subset of all choices $(X_1,\ldots,X_{p-1})$ in $\Omega$ such that there is a small pair containing $u$ and a vertex in $\cup_{i\in[p-1]}X_i$.

\begin{claim}\label{runiformclaimnew} $\Omega=\Omega'\cup (\cup_{i\in[4]}\Omega_{u_i})$.
\end{claim}
\begin{prof}
Let   $\mathcal{X}=(X_1,\ldots,X_{p-1})\in \Omega$, where $X_i=\{x_{i1},\ldots,x_{i,v(F)},y_{i1},\ldots,y_{i,r-2}\}$ for $i\in [2]$ and $X_i=\{x_{i1},\ldots,x_{i,v(F)}\}$ for $i\in [p-1]\backslash \{1,2\}$. We split the proof into several cases depending on $\ell$ and the positions of $\{u_1,u_2\},\{u_3,u_4\}$.

\noindent{\bf Case 1.} \quad $\{u_1,u_2\},$ $\{u_3,u_4\}$  are in different vertex classes.

Without loss of generality, assume $\{u_1,u_2\}\subseteq V_1,$ $\{u_3,u_4\}\subseteq V_2$.

\noindent{\bf Subcase 1.1}\quad $\ell=2$.

 Consider the two edges $e_1'=\{x_{11},x_{12},y_{11},\ldots,y_{1,r-2}\}$, $e_2'=\{x_{21},x_{22},y_{21},\ldots,y_{2,r-2}\}$ in $K_n^r$. If $C(e_1')=C(e_2)$ and $C(e_2')=C(e_1)$, then $C(e_1')\neq C(e_2')$. Thus, there is a small pair in ${{\cup_{i\in[p-1]}X_i}\choose 2}$, that is, $\mathcal{X}\in \Omega'$. Otherwise, let $X_i'=\{x_{i1},\ldots,x_{i,q_i}\}$ for $i\in[p-1]$. Then, $\cup_{i\in[p-1]}X_i'$  could form the skeleton of $F^{(r)}$, and then be expanded to a rainbow copy of $F^{(r)}$ containing $e_1'$ and $e_2'$, a contradiction. If $C(e_1')\neq C(e_2)$, then  there is a small pair consisting of a vertex in $\{u_3,u_4\}$ and another vertex in $\cup_{i\in[p-1]}X_i$, that is, $\mathcal{X}\in \Omega_{u_3}\cup \Omega_{u_4}$. Otherwise, let $X_i'=\{x_{i1},\ldots,x_{i,q_i}\}$ for $i\in [p-1]\backslash \{2\}$ and $X_2'=\{x_{21},\ldots,x_{2,q_2-2}\}$. Then, $\cup_{i\in[p-1]}X_i'\cup\{u_3,u_4\}$ could form the skeleton of $F^{(r)}$, and then be expanded to a rainbow copy of $F^{(r)}$ containing $e_1'$ and $e_2$, a contradiction. The case $C(e_2')\neq C(e_1)$ is similar.

Thus, $\Omega=\Omega'\cup (\cup_{i\in[4]}\Omega_{u_i})$.

\noindent{\bf Subcase 1.2}\quad $ \ell\geq 3$ and the two edges of $F$ in $U_1'$ intersect.

 Consider the  edge $e=\{u_1,x_{11},u_3,x_{21},y_{11}\ldots,y_{1,r-4}\}$  in $K_n^r$. Then, at least one of $\{e,e_1\} $ and $\{e,e_2\}$ is rainbow. Moveover, both $\{e,e_1\} $ and $\{e,e_2\}$ are pairs of non-crossing edges which contain intersecting non-crossing parts in a common vertex class. Thus, $\mathcal{X}\in \Omega'\cup (\cup_{i\in[4]}\Omega_{u_i})$. Otherwise, let $X_i'=\{x_{i1},\ldots,x_{i,q'_i}\}$ for $i\in[p-1]$ and $X_i''=\{x_{i1},\ldots,x_{i,q'_i-2}\}$ for $i=1,2$.   It follows that both $X_1''\cup(\cup_{i\in[p-1]\backslash \{1\}}X_i')\cup \{u_1,u_2\}$ and $X_2''\cup(\cup_{i\in[p-1]\backslash \{2\}}X_i')\cup \{u_3,u_4\}$ could form the skeleton of $F^{(r)}$. Thus, at least one of them can be expanded to a rainbow copy of $F^{(r)}$ containing either $\{e, e_1\}$ or $\{e,e_2\},$   a contradiction.

\noindent{\bf Subcase 1.3}\quad $\ell\geq 3$ and the two edges of $F$ in $U_1'$ are disjoint.

 Consider the  edge $e=\{x_{11},x_{12},x_{21},x_{22},y_{11}\ldots,y_{1,r-4}\}$  in $K_n^r$. Then, at least one of $\{e,e_1\} $ and $\{e,e_2\}$ is rainbow.  Thus, $\mathcal{X}\in \Omega'\cup (\cup_{i\in[4]}\Omega_{u_i})$. Otherwise, let $X_i'=\{x_{i1},\ldots,x_{i,q'_i}\}$ for $i\in[p-1]$ and $X_i''=\{x_{i1},\ldots,x_{i,q'_i-2}\}$ for $i=1,2$.   It follows that both $X_1''\cup(\cup_{i\in[p-1]\backslash \{1\}}X_i')\cup \{u_1,u_2\}$ and $X_2''\cup(\cup_{i\in[p-1]\backslash \{2\}}X_i')\cup \{u_3,u_4\}$ could form the skeleton of $F^{(r)}$. Thus at least one of them can be expanded to a rainbow copy of $F^{(r)}$ containing either $\{e, e_1\}$ or $\{e,e_2\}$, a contradiction.

\noindent{\bf Case 2.} \quad $\{u_1,u_2\},$ $\{u_3,u_4\}$ are contained in a common vertex class.

Without loss of generality, assume $\{u_1,\ldots,u_{4}\} \subseteq V_1$ and $u_1\neq u_3$, $u_1\neq u_4$, $u_2\neq u_3$.

\noindent{\bf Subcase 2.1} \quad $\ell= 2$.

Consider $e=\{x_{21},x_{22},y_{21},\ldots,y_{2,r-2}\}$. At least one of $\{e,e_1\}$ and  $\{e,e_2\}$ is rainbow. Thus, $\mathcal{X}\in\Omega'\cup (\cup_{i\in[4]}\Omega_{u_i})$. Otherwise, let $X_i'=\{x_{i1},\ldots,x_{i,q_i}\}$ for $i\in [p-1]\backslash \{1\}$ and $X_1'=\{x_{11},\ldots,x_{1,q_1-2}\}$. It follows  that   both  $\cup_{i\in[p-1]}X_i'\cup \{u_1,u_2\}$ and $\cup_{i\in[p-1]}X_i'\cup \{u_3,u_4\}$ could form the skeleton of $F^{(r)}$. Thus at least one of them can be expanded to a rainbow copy of $F^{(r)}$ containing either $\{e, e_1\}$ or $\{e,e_2\},$  a contradiction.

\noindent{\bf Subcase 2.2} \quad $\ell\geq 3$ and the two edges of $F$ in $U_1'$ intersect.

Consider $e=\{u_1,u_3,y_{11},\ldots,y_{1,r-2}\}$. Then at least one of $\{e,e_1\}$ and  $\{e,e_2\}$ is rainbow. Thus, $\mathcal{X}\in\Omega'\cup (\cup_{i\in[4]}\Omega_{u_i})$. Otherwise, let $X_i'=\{x_{i1},\ldots,x_{i,q'_i}\}$ for $i\in [p-1]\backslash \{1\}$ and $X_1'=\{x_{11},\ldots,x_{1,q'_1-3}\}$. It follows  that both  $\cup_{i\in[p-1]}X_i'\cup \{u_1,u_2,u_3\}$ and $\cup_{i\in[p-1]}X_i'\cup \{u_1,u_3,u_4\}$ could form the skeleton of $F^{(r)}$. Thus at least one of them can be expanded to a rainbow copy of $F^{(r)}$ containing either $\{e, e_1\}$ or $\{e,e_2\},$  a contradiction.

\noindent{\bf Subcase 2.3}\quad $\ell\geq 3$ and the two edges of $F$ in $U_1'$ are disjoint.

  Consider $e=\{x_{11},x_{12},y_{11},\ldots,y_{1,r-2}\}$. Then at least one of $\{e,e_1\}$ and  $\{e,e_2\}$ is rainbow. Thus, $\mathcal{X}\in\Omega'\cup (\cup_{i\in[4]}\Omega_{u_i})$. Otherwise, let $X_i'=\{x_{i1},\ldots,x_{i,q_i}\}$ for $i\in [p-1]\backslash \{1\}$ and $X_1'=\{x_{11},\ldots,x_{1,q_1-2}\}$. It follows  that   both  $\cup_{i\in[p-1]}X_i'\cup\{u_1,u_2\}$ and $\cup_{i\in[p-1]}X_i'\cup\{u_3,u_4\}$ could form the skeleton of $F^{(r)}$. Thus at least one of them can be expanded to a rainbow copy of $F^{(r)}$ containing either $\{e, e_1\}$ or $\{e,e_2\},$  a contradiction.

Therefore, in all cases, we have $\Omega=\Omega'\cup (\cup_{i\in[4]}\Omega_{u_i})$.
\end{prof}

For $v\in V(G)$, define $d_s(v)=|\{u\in V(G): \{u,v\} \mbox{ is small}\}|.$ Let $A=\{v\in V(G):d_s(v)\geq \varepsilon_3 n\}$.

\begin{claim}\label{runiformclaim2} Let $\{e_1,e_2\} $ be any pair of non-crossing edges, and $\{u_1,u_2\},$ $\{u_3,u_4\}$ be a   non-crossing part of $e_1,e_2$, respectively. We have at least one of the vertices $u_i$, $i\in[4]$ is in $A$.
\end{claim}
\begin{prof}
  By Claim \ref{runiformclaimnew}, $\Omega=\Omega'\cup (\cup_{i\in[4]}\Omega_{u_i})$.
Then, we bound the size of $\Omega_{u_i}$. Clearly, $$|\Omega|\geq   \bigg({{\frac{n}{3(p-1)}}\choose {v(F)+r-2}}\bigg)^2 \bigg({{\frac{n}{3(p-1)}}\choose {v(F)}}\bigg)^{p-3}\geq \varepsilon_2 n^{(p-1)v(F)+2r-4}.$$
Since each small pair contained in $\Omega'$ corresponds to at most $n^{(p-1)v(F)+2r-6}$ choices in $\Omega'$, combined with Claim \ref{runiformclaim1}, we have
$$|\Omega'|\leq \varepsilon_6 n^2\cdot n^{(p-1)v(F)+2r-6}\leq \frac{1}{2}|\Omega|.$$
Therefore,
$$\sum_{i\in[4]}|\Omega_{u_i}|\geq \frac{1}{2}|\Omega|,$$
which follows that there is some vertex $u\in\{u_i:i\in [4]\}$ such that $|\Omega_u|\geq \frac{1}{8}|\Omega|\geq \frac{\varepsilon_2}{8}n^{(p-1)v(F)+2r-4}$. Hence,
$$d_s(u)\geq \frac{\varepsilon_2 n^{(p-1)v(F)+2r-4}}{8 n^{(p-1)v(F)+2r-5}}\geq \varepsilon_3 n.$$

Therefore,  at least one of the vertices $u_i$, $i\in[2\ell]$ is in $A$.
\end{prof}

Now, let $G'$ be the subgraph of $G$ obtained by deleting from $G$ all the non-crossing edges whose non-crossing parts are disjoint from $A$. By Claim \ref{runiformclaim2}, the number of deleted edges is at most 1.

For any $v\in A$, define $N_{ncp}(v)=\{u\in V(G'): \{u,v\} \mbox{ is a non-crossing part of some non-crossing edge in } G'\}$ and $d_{ncp}(v)=|N_{ncp}(v)|$.

\begin{claim}\label{runifclaim3}
There is a vertex $v_0\in A$ such that the number of non-crossing edges whose non-crossing part containing $v_0$ is at least  $\varepsilon_4 n^{r-1}$ and $d_{ncp}(v_0)\geq \varepsilon_5 n$.
\end{claim}
\begin{prof} By the definition of $A$, the number of small pairs intersecting $A$ is at least $|A|\varepsilon_3 n/2$. Then,
$$|T\backslash G|\geq \frac{|A|\varepsilon_3 n}{2}\bigg( {{p-3}\choose {r-2}}\bigg(\frac{n}{2(p-1)}\bigg)^{r-2}-\bigg({n\choose {r-3}}\tau+|F|\bigg)  \bigg)\frac{1}{{r\choose 2}}\geq \varepsilon_4 |A| n^{r-1}.$$
Thus,
$$|G'\backslash T|\geq |G\backslash T|-1\geq |T\backslash G|+2-1\geq \varepsilon_4 |A| n^{r-1}.$$
Since the non-crossing part of every edge in $G'\backslash T$ intersects $A$,  there is a vertex $v_0\in A$ such that the number of non-crossing edges whose non-crossing part containing $v_0$ is at least  $\varepsilon_4 n^{r-1}$. And this implies $d_{ncp}(v_0)\geq \varepsilon_4 n^{r-1}/n^{r-2}\geq \varepsilon_5 n$.
\end{prof}

Without loss of generality, assume $v_0\in V_1$. For $2\leq j\leq p-1$, let $Z_j=\{v\in V_j:\{v_0,v\} \mbox{ is big}\}$. For $y\in N_{ncp}(v_0)$, let $e_{v_0}(\cdot)$ be any mapping from $N_{ncp}(v_0)$ to the set of non-crossing edges containing $v_0$ such that $e_{v_0}(y)$ contains $y$. Now we finish the proof by discussion on the size of each $Z_j$, $2\leq j \leq p-1$.

\noindent{\bf Case A.} \quad $|Z_j|\geq \varepsilon_5 n$ for all $2\leq j \leq p-1$.

\noindent{\bf Subcase A.1}\quad  $\ell\geq 3$ and the two edges of $F$ in $U_1'$ intersect.

Let \begin{align*}&\Omega=\bigg\{(y,z,W,W_1,\ldots,W_{p-1})\in N_{ncp}(v_0)\times V_1\times {V\choose {r-2}}\times{V_1\choose {q'_1-3}}\times{Z_2\choose {q'_2}}\times\cdots \times{Z_{p-1}\choose {q'_{p-1}}}:\\
& \quad\quad \ \ \  (W\cup \{z\})\cap  e_{v_0}(y)=\emptyset, (W\cup \{z\}\cup e_{v_0}(y))\cap (\cup_{i\in[p-1]}W_i)=\emptyset \bigg\},\\
&\Omega_0=\{(y,z,W,W_1,\ldots,W_{p-1})\in \Omega:C(e_{v_0}(y))=C(\{v_0,z,W\})\},\\
& \Omega_1=\Omega\backslash\Omega_0.
\end{align*}

First, $$|\Omega|\geq\varepsilon_5 n\frac{n}{3(p-1)}{{n/2}\choose {r-2}}{{\frac{n}{3(p-1)}}\choose {q'_1-3}}\Pi_{i=2}^{p-1}{{\frac{\varepsilon_5n}{2}}\choose {q'_i}}>2\varepsilon_6n^{v(F)+r-3}.$$
Next, we bound $|\Omega_0|.$ Recall that $G$ is rainbow. Thus, for a fixed choice of $\{z,W\}$, there are at most $(r-1)$ choices of $y$ such that $C(e_{v_0}(y))=C(\{v_0,z,W\})$. Therefore, $|\Omega_0|\leq (r-1){n\choose {r-1}} n^{v(F)-3}\leq \frac{|\Omega|}{2}$.
As $\Omega=\Omega_0\cup\Omega_1$, we have $|\Omega_1|\geq |\Omega|/2$. Since $C(e_{v_0}(y))\neq C(\{v_0,z,W\})$ holds for every choice in $\Omega_1$, there is at least one small pair which does not contain $v_0$  in each choice in $\Omega_1$. Each such small pair is counted at most $n^{v(F)+r-5}$ times in $\Omega_1$. Thus, the total number of small pairs is at least
$$\frac{|\Omega_1|}{n^{v(F)+r-5}}\geq \frac{|\Omega|}{2n^{v(F)+r-5}}> \varepsilon_6 n^2 ,$$
which contradicts to Claim \ref{runiformclaim1}.

\noindent{\bf Subcase A.2}\quad  $\ell\geq 3$ and the two edges of $F$ in $U_1'$ are disjoint.

Similar to Subcase A.1, let \begin{align*}&\Omega=\bigg\{(y,Z,W,W_1,\ldots,W_{p-1})\in N_{ncp}(v_0)\times {V_1\choose 2}\times {V\choose {r-2}}\times{V_1\choose {q'_1-4}}\times{Z_2\choose {q'_2}}\times\cdots \times{Z_{p-1}\choose {q'_{p-1}}}:\\
&\quad\quad \ \ \  (W\cup Z)\cap  e_{v_0}(y)=\emptyset, (W\cup Z\cup e_{v_0}(y))\cap (\cup_{i\in[p-1]}W_i)=\emptyset \bigg\},\\
&\Omega_0=\{(y,Z,W,W_1,\ldots,W_{p-1})\in \Omega:
C(e_{v_0}(y))=C(\{Z,W\})\},\\
& \Omega_1=\Omega\backslash\Omega_0.
\end{align*}

First, $$|\Omega|\geq\varepsilon_5 n{{\frac{n}{3(p-1)}}\choose 2}{{n/2}\choose {r-2}}{{\frac{n}{3(p-1)}}\choose {q'_1-4}}\Pi_{i=2}^{p-1}{{\frac{\varepsilon_5n}{2}}\choose {q'_i}}>2\varepsilon_6n^{v(F)+r-3}.$$
Next, we bound $|\Omega_0|.$ Recall that $G$ is rainbow. Thus, for a fixed choice of $\{Z,W\}$, there are at most $(r-1)$ choices of $y$ such that $C(e_{v_0}(y))=C(\{Z,W\})$. Therefore, $|\Omega_0|\leq (r-1){n\choose {r}} n^{v(F)-4}\leq \frac{|\Omega|}{2}$.
As $\Omega=\Omega_0\cup\Omega_1$, we have $|\Omega_1|\geq |\Omega|/2$. Since $C(e_{v_0}(y))\neq C(\{Z,W\})$ holds for every choice in $\Omega_1$, there is at least one small pair which does not contain $v_0$  in each choice in $\Omega_1$. Each such small pair is counted at most $n^{v(F)+r-5}$ times in $\Omega_1$. Thus, the total number of small pairs is at least
$$\frac{|\Omega_1|}{n^{v(F)+r-5}}\geq \frac{|\Omega|}{2n^{v(F)+r-5}}> \varepsilon_6 n^2 ,$$
which contradicts to Claim \ref{runiformclaim1}.

\noindent{\bf Subcase A.3}\quad  $\ell=2$.

In this case, we use the partition $\mathcal{Q}$ of $F$.
Similarly, let \begin{align*}&\Omega=\bigg\{(y,Z,W,W_1,\ldots,W_{p-1})\in N_{ncp}(v_0)\times {Z_2\choose 2}\times {V\choose {r-2}}\times{V_1\choose {q_1-2}}\times{Z_2\choose {q_2-2}}\times {Z_3\choose {q_3}}\times\cdots \times{Z_{p-1}\choose {q_{p-1}}}:\\
&\quad\quad \ \ \  (W\cup Z) \cap  e_{v_0}(y)=\emptyset, (W\cup Z \cup e_{v_0}(y))\cap (\cup_{i\in[p-1]}W_i)=\emptyset \bigg\},\\
&\Omega_0=\{(y,Z,W,W_1,\ldots,W_{p-1})\in \Omega:
C(e_{v_0}(y))=C(\{Z,W\})\},\\
& \Omega_1=\Omega\backslash\Omega_0.
\end{align*}

First, $$|\Omega|\geq\varepsilon_5 n{{\frac{\varepsilon_5n}{2}}\choose {2}}{{n/2}\choose { r-2}}{{\frac{n}{3(p-1)}}\choose {q_1-2}}{{\frac{\varepsilon_5n}{2}}\choose {q_2-2}}\Pi_{i=3}^{p-1}{{\frac{\varepsilon_5n}{2}}\choose {q_i}}>2\varepsilon_6n^{v(F)+r-3}.$$
Next, we bound $|\Omega_0|.$ Again, notice that $G$ is rainbow. Thus, for a fixed choice of $\{Z,W\}$, there are at most $(r-1)$ choices of $y$ such that $C(e_{v_0}(y))=C(\{Z,W\})$. Therefore, $|\Omega_0|\leq (r-1){n\choose {r}} n^{v(F)-4}  \leq \frac{|\Omega|}{2}$.
As $\Omega=\Omega_0\cup\Omega_1$, we have $|\Omega_1|\geq |\Omega|/2$. Since $C(e_{v_0}(y))\neq C(\{Z,W\})$ holds for every choice in $\Omega_1$, there is at least one small pair which does not contain $v_0$  in each choice in $\Omega_1$. Each such small pair is counted at most $n^{v(F)+r-5}$ times in $\Omega_1$. Thus, the total number of small pairs is at least
$$\frac{|\Omega_1|}{n^{v(F)+r-5}}\geq \frac{|\Omega|}{2n^{v(F)+r-5}}> \varepsilon_6 n^2 ,$$
which contradicts to Claim \ref{runiformclaim1}.

\noindent{\bf Case B.} \quad There is some $j\in[p-1]\backslash \{1\}$ such that $|Z_j|< \varepsilon_5 n$.

Without loss of generality, assume $j=2$.
We show that if $v_0$ is moved from $V_1$ to $V_2$, then $f_G(\mathcal{P^*})$ will strictly increase, which contradicts to that $f_G(\mathcal{P^*})$ is maximum. Let $E'=\{e\in G: e\cap V_1=\{v_0\},e\cap V_2\neq \emptyset\}$, $E''=\{e\in G: v_0\in e, |e\cap V_1|\geq 2,e\cap V_2= \emptyset\}$. It suffices to show $|E''|>|E'|$.
As $|Z_2|< \varepsilon_5 n$, we have
$$|E'|\leq \varepsilon_5 n{n\choose {r-2}}+(n-\varepsilon_5 n)\bigg({n\choose {r-3}}\tau+|F|\bigg)\leq 2\varepsilon_5 n^{r-1}.$$
To bound $|E''|$, recall that the number of non-crossing edges with non-crossing part containing $v_0$ is at least $\varepsilon_4 n^{r-1}$ by Claim \ref{runifclaim3}. Moreover, the number of edges containing $v_0$ and a vertex in $V_2$ is at most $  \varepsilon_5 n{n\choose {r-2}}+(n-\varepsilon_5 n)({n\choose {r-3}}\tau+|F|)\leq 2\varepsilon_5 n^{r-1}.$ Thus, $|E''|\geq \varepsilon_4 n^{r-1}-2\varepsilon_5 n^{r-1}>2\varepsilon_5 n^{r-1}=|E'|$. The proof is completed.
\end{prof}

\section*{Reference}

\bibliography{bibl}

\begin{thebibliography}{10}
\expandafter\ifx\csname urlstyle\endcsname\relax
  \providecommand{\doi}[1]{doi:\discretionary{}{}{}#1}\else
  \providecommand{\doi}{doi:\discretionary{}{}{}\begingroup
  \urlstyle{rm}\Url}\fi

\bibitem{Erdos&Simonovits&Sos1975}
P.~Erd\H{o}s, M.~Simonovits and V.~T. S$\acute{\rm o}$s.
\newblock \emph{Anti-{Ramsey} theorems}.
\newblock in Infinite and FiniteSets, Vol. II. Colloq. Math. Soc. J\'{a}nos
  Bolyai, \textbf{10 (1975)}, 633--643.

\bibitem{Frankl&Kupavskii2019}
P.~Frankl and A.~Kupavskii.
\newblock \emph{Two problems of {P. Erd$\mbox{\H{o}}$s} on matchings in set
  families---in the footsteps of {Erd$\mbox{\H{o}}$s and Kleitman}}.
\newblock J. Combin. Theory Ser. B, \textbf{138 (2019)}, 286--313.

\bibitem{Fujita&Magnant&Ozeki2010}
S.~Fujita, C.~Magnant and K.~Ozeki.
\newblock \emph{Rainbow generalizations of {Ramsey} theory: A survey}.
\newblock Graphs Combin., \textbf{26 (2010)}, 1--30.

\bibitem{Gu&Li&Shi2020}
R.~Gu, J.~Li and Y.~Shi.
\newblock \emph{Anti-ramsey numbers of paths and cycles in hypergraphs}.
\newblock SIAM Journal on Discrete Mathematics, \textbf{34 (2020), no.~1},
  271--307.

\bibitem{Jiang&Pikhurko2009}
T.~Jiang and O.~Pikhurko.
\newblock \emph{Anti-ramsey numbers of doubly edge-critical graphs}.
\newblock Journal of Graph Theory, \textbf{61 (2009), no.~3}, 210--218.

\bibitem{Montellano-Ballesteros&Neumann-Lara2002}
J.~J. Montellano-Ballesteros and V.~Neumann-Lara.
\newblock \emph{An anti-{Ramsey} theorem}.
\newblock Combinatorica, \textbf{22 (2002)}, 445--449.

\bibitem{Mubayi1}
D.~Mubayi.
\newblock \emph{A hypergraph extension of {T}ur\'{a}n's theorem}.
\newblock J. Combin. Theory Ser. B, \textbf{96 (2006)}, 122--134.

\bibitem{Ozkahya&Young2013}
L.~\"{O}zkahya and M.~Young.
\newblock \emph{Anti-{Ramsey} number of matchings in hypergraphs}.
\newblock Discrete Math., \textbf{313 (2013)}, 2359--2364.

\bibitem{Pikhurko}
O.~Pikhurko.
\newblock \emph{Exact computation of the hypergraph tur$\acute{\rm a}$n
  function for expanded complete 2-graphs}.
\newblock J. Combin. Theory Ser. B, \textbf{103 (2013)}, 220--225.

\bibitem{Tang&Li&Yan2022}
Y.~Tang, T.~Li and G.~Yan.
\newblock \emph{Anti-ramsey number of expansions of paths and cycles in uniform
  hypergraphs}.
\newblock Journal of Graph Theory, \textbf{101 (2022), no.~4}, 668--685.

\bibitem{Tang&Li&Yan2022+}
Y.~C. Tang, T.~Li and G.~Y. Yan.
\newblock \emph{Anti-ramsey numbers of cycles of length three in uniform
  hypergraphs}, 2022.
\newblock In press.

\end{thebibliography}

\end{document}